\documentclass[12pt]{amsart}

\usepackage{mathrsfs}

\usepackage{amsmath,amsthm,amssymb}
\font\teneufm=eufm10
\font\seveneufm=eufm7
\font\fiveeufm=eufm5
\newfam\frakturfam

\textfont\frakturfam=\teneufm
\scriptfont\frakturfam=\seveneufm
\scriptscriptfont\frakturfam=\fiveeufm

\newtheorem{pr}{Proposition}

\newtheorem{lm}{Lemma}
\newtheorem{theor}{Theorem}
\newtheorem{co}{Corollary}

\def\bee{\begin{eqnarray}}
\def\bes{\begin{eqnarray*}}
\def\eee{\end{eqnarray}}
\def\ees{\end{eqnarray*}}

\def\Proof{{\sl Proof.}\ }

\title{Differential algebraic dependence and Novikov dependence}

\begin{document}
\date{}
\maketitle

\begin{center}
{\bf Bibinur Duisengaliyeva,}\footnote{
 L.N. Gumilyov Eurasian National University,
Astana, 010008, Kazakhstan,
e-mail: {\em bibinur.88@mail.ru}}
{\bf Ualbai Umirbaev}\footnote{
 Wayne State University,
Detroit, MI 48202, USA, and Institute of Mathematics and Modeling, Almaty, Kazakhstan, 
e-mail: {\em umirbaev@math.wayne.edu}

The project is partially supported by the grant AP05133009 of MES RK}

\end{center}

\begin{abstract}  We define an analogue of the Fox derivatives for differential polynomial algebras and give a criterion for differential algebraic dependence of a finite system of elements. In particular, we prove that differential algebraic dependence of a finite set of elements of a differential polynomial algebra over a constructive differential field $k$ of characteristic zero is algorithmically recognizable. Using a representation of free Novikov algebras by differential polynomials we also give a criterion of Novikov dependence of a finite system of elements of free Novikov algebras.
\end{abstract}

\noindent


\section{Introduction}

\hspace*{\parindent}

Algebraic dependence of a finite set of elements of a polynomial algebra $k[x_1,x_2,\ldots,x_n]$ over a constructive field $k$ is algorithmically recognizable  \cite{SS} by the methods of Gr\"obner bases \cite{CLO}. The concept of algebraic dependence is easily generalized to arbitrary varieties of algebras. For example, a system of elements $f_1,f_2,\ldots,f_m$ of an arbitrary associative algebra $A$ is called {\em associatively independent} if the subalgebra generated by these elements is a free associative algebra with free generators $f_1,f_2,\ldots,f_m$. Otherwise the elements $f_1,f_2,\ldots,f_m$ are called {\em associatively dependent}.

By the Nielsen-Schreier Theorem the subgroups of free groups are free \cite{MKS} and by the Shirshov-Witt Theorem the subalgebras of free Lie algebras are free \cite{Shir1, Witt}. These results easily imply that the dependence of a finite system of elements of free groups and free Lie algebras are algorithmically  recognizable. Associative dependence of a finite system of elements of free associative algebras \cite{93AL2} is algorithmically unrecognizable.

This paper is devoted to the study of differential algebraic dependence and Novikov dependence. The basic concepts of differential algebras can be found in \cite{Kolchin,Ritt}. Let $\Delta=\{\delta_1,\ldots,\delta_m\}$ be a basic set of derivation operators. If $\Delta$ contains only one element, then differential algebras are called {\em ordinary}, and if $\Delta$ contains at least two elements, then they are called {\em partial}.

The differential polynomial algebra $k\{x_1,x_2,\ldots,x_n\}$ in the variables $x_1,x_2,\ldots,x_n$ over a differential field $k$ is the closest generalization of the polynomial algebra $k[x_1,x_2,\ldots,x_n]$. However, the structures of the subalgebras and ideals of these algebras are very different. It is well known that the ideal membership problem (see, for example \cite{CLO}) and the subalgebra membership problem \cite{SS,Noskov} for polynomial algebras are algorithmically decidable. At the same time the ideal membership problem and the subalgebra membership problem for partial differential polynomial algebras is algorithmically undecidable \cite{16JA}. These questions remain open for ordinary differential polynomial algebras \cite{KZ}.

A set of elements $f_1,f_2,\ldots,f_m$ of the polynomial algebra $k[x_1,\ldots,x_n]$ over a field $k$ of characteristic zero is algebraically dependent if and only if the rank of the Jacobian matrix
\bes
J(f_1,f_2,\ldots,f_m)=[\frac{\partial f_i}{\partial x_j}]_{1\leq i\leq m, 1\leq j\leq n}
\ees
is less than $m$ (see, for example \cite{SU4}). In this paper we formulate and prove an analogue of this result for the differential polynomial algebra $k\{x_1,x_2,\ldots,x_n\}$ over a differential field $k$ of characteristic zero in the terms of Fox derivatives (see, for example \cite{12JA}). This result allows us to prove that the differential algebraic dependence of a finite system of elements of differential polynomial algebras over a constructive field of characteristic zero is algorithmically recognizable.

I.M. Gel'fand and I.Ya. Dorfman \cite{13} noticed that any differential algebra with derivation $\delta$ with respect to the multiplication $a\circ b=a(\delta b)$ becomes a Novikov algebra. In \cite{DzhL} free Novikov algebras are represented by ordinary differential polynomial algebras via the multiplication $a\circ b=a(\delta b)$. Using this representation of free Novikov algebras, we show that a finite system of elements of a free Novikov algebra is Novikov dependent if and only if it is differentially algebraically dependent. In particular, the Novikov dependence of a finite system of elements of free Novikov algebras is also algorithmically recognizable. This also implies that any $n+1$ elements of a free Novikov algebra of rank $n$ are Novikov dependent.

It is well known \cite{Essen} that any two algebraically dependent elements of a polynomial algebra are polynomials in one variable. Any two associatively  dependent elements of a free associative algebra are also polynomials in one variable \cite{Cohn}. L. Makar-Limanov and I. Shestakov \cite{MLSh} (see also \cite{12JA,MLU16}) proved an analogue of this result for free Poisson algebras over fields of characteristic zero. The question on the validity of similar results for differential polynomial algebras and for free Novikov algebras remains open.

The rest of the paper is organized as follows. In Section 2 we give some elementary terminology on differential polynomial algebras. The definition of the Fox derivatives for differential polynomial algebras is given in Section 3. Section 3 is devoted to the study of differential algebraic dependence in differential polynomial algebras. In Section 4 we recall a representation of free Novikov algebras in differential polynomial algebras. Using this reperesentation, in Section 5, we describe Novikov dependence of elements in free Novikov algebras.

\section{Differential polynomial algebras}

\hspace*{\parindent}

Let $\Delta = \{\delta_1,\ldots, \delta_m\}$ be a basic set of derivation operators. A commutative ring $R$ with identity is called a {\em differential ring} or {\em $\Delta$-ring} if all elements of $\Delta$ act on $R$ as a commuting set of derivations, i.e., the derivations $\delta_i: R\rightarrow R$ are defined for all $i$ and $\delta_i\delta_j=\delta_j\delta_i$ for all $i,j$.

Let $\Theta$ be the free commutative monoid on the set $\Delta = \{\delta_1,\ldots, \delta_m\}$  of derivation operators. The elements
\bes
\theta = \delta_1^{i_1}\ldots \delta_m^{i_m}
\ees
of the monoid $\Theta$ are called {\em derivative operators}. The {\em order} of $\theta$ is defined as $|\theta|=i_1+\ldots+i_m$. Let
$\gamma(\theta)=(i_1,\ldots,i_m)\in \mathbb{Z}_+^m$, where $\mathbb{Z}_+$ is the set of all nonnegative integers.

Let $R$ be a differential ring. Denote by $R[\Delta]$ the free left $R$-module with a basis $\Theta$. Every element $u\in R[\Delta]$ can be uniquely written in the form
\bes
u=\sum_{\theta\in \Theta} r_{\theta}\theta
\ees
with a finite number of nonzero $r_{\theta}\in R$. We turn $R[\Delta]$ to a ring by
\bes
\delta_i r=r \delta_i+ \delta_i(r)
\ees
for all $i$ and $r\in R$. It is well known \cite{KLMP} that these relations uniquely define a structure of a ring on $R[\Delta]$. The ring $R[\Delta]$ is called the ring of {\em differential operators with coefficients in $R$}. Every left module over $R[\Delta]$ is called a {\em differential module} over $R$. In particular, $R$ is a left $R[\Delta]$ and every $I\subseteq R$ is a differential ideal of $R$ if and only if $I$ is an $R[\Delta]$-submodule of $R$. This means that $R[\Delta]$ is the universal enveloping ring of $R$ in the usual ring-theoretic terminology.

Let $x^\Theta=\{x^\theta | \theta\in \Theta\}$ be a set of symbols enumerated by the elements of $\Theta$. Consider the polynomial algebra $R[x^\Theta]$ over $R$ generated by the set of (polynomially) independent variables $x^\Theta$. It is easy to check that the derivations $\delta_i$ can be uniquely extended to a derivation of $R[x^\Theta]$ by $\delta_i(x^\theta)=x^{\delta_i\theta}$. Denote this differential ring
by $R\{x\}$; it is called the {\em ring of differential polynomials} in $x$ over $R$.

By adjoining more variables, we can obtain the differential ring $R\{x_1,x_2,\ldots,x_n\}$ of the differential polynomials in $x_1,x_2,\ldots,x_n$ over $R$. Let $M$ be the free commutative monoid generated by all elements $x_i^{\theta}$, where $1\leq i\leq n$ and $\theta\in \Theta$. The elements of $M$ are called {\em differential monomials} of $R\{x_1,x_2,\ldots,x_n\}$. Every element $a\in R\{x_1,x_2,\ldots,x_n\}$ can be uniquely written in the form
\bes
a=\sum_{m\in M} r_m m
\ees
with a finite number of nonzero $r_m\in R$.

If $R$ is a domain then the field of fractions $R\langle x_1,x_2,\ldots,x_n\rangle$ of $R\{x_1,x_2,\ldots,x_n\}$ becomes a differential field since every $\delta_i$ can be uniquely extended to $R\langle x_1,x_2,\ldots,x_n\rangle$. This field is called the field of {\em rational differential functions} in the variables $x_1,x_2,\ldots,x_n$ over $R$.

\section{Fox derivatives}

\hspace*{\parindent}

Let $k$ be an arbitrary differential field of characteristic zero and let $A=k\{x_1,x_2,\ldots,x_n\}$ be the differential polynomial algebra over $k$ in the variables $x_1,x_2,\ldots,x_n$. Let $B=k\langle x_1,x_2,\ldots,x_n\rangle$. We are going to define analogues of the Fox derivatives for $A$ and $B$. Let $C=k\langle x_1,x_2,\ldots,x_n,y_1,y_2,\ldots,y_n\rangle$ be the differential field with additional variables $y_1,y_2,\ldots,y_n$ and $B\subset C$. Let $\lambda$ be the derivation of the differential field $C$ defined by $x_i\mapsto y_i$ and $y_i\mapsto 0$ for all $i$. Put
\bes
\Omega_B=B[\Delta]y_1\oplus B[\Delta]y_2\oplus\ldots\oplus B[\Delta]y_n,
\ees

Notice that if $f\in B$ then $\lambda(f)\in \Omega_B$. Denote by
\bes
D : B\longrightarrow \Omega_B
\ees
the restriction of $\lambda$.  Then $D$ is a derivation of $B$ with coefficients in the free left $B[\Delta]$-module $\Omega_B$. Obviously, $D$ is the universal derivation of $B$ (see, for example \cite{12JA}). For any $f\in B$ there are unique elements $a_1,\ldots,a_n\in B[\Delta]$ such that
\bes
D(f)=a_1y_1+\ldots+a_ny_n.
\ees
Let $a_i=d_{x_i}(f)$ for all $i$. This defines the derivations
\bes
d_{x_i} : B\longrightarrow B[\Delta]
\ees
such that $d_{x_i}(x_j)=\delta_{ij}$, where $\delta_{ij}$ is the Kronecker delta. These derivations are analogues of Fox derivatives. Let
\bes
d(f)=(d_{x_1}(f),\ldots,d_{x_n}(f)), \ \ f\in B.
\ees
If $f_1,\ldots,f_m\in B$ then the matrix
\bes
J(f_1,\ldots,f_m)=[d_{x_j}(f_i)]_{1\leq i\leq m, 1\leq j\leq n}
\ees
is called the {\em Jacobian matrix} of  $f_1,\ldots,f_m$. The row vectors $d(f_1),\ldots,d(f_m)$ are rows of $J(f_1,\ldots,f_m)$. The Chain Rule can be written in the form
\bes
d(f(g_1,\ldots,g_n))=d(f)J(g_1,\ldots,g_n).
\ees

Notice that if $f\in A$ then $d_{x_i}(f)\in A[\Delta]$. Let $\varphi$ be an endomorphism of $A$ that sends $x_i$ to $f_i$. Then
\bes
J(\varphi)=J(f_1,\ldots,f_n)
\ees
is the {\em Jacobian matrix} of $\varphi$. The Chain Rule implies that $J(\varphi)$ is invertible over $A[\Delta]$ if $\varphi$ is an automorphism. The Jacobian Conjecture for differential polynomial algebras can be formulated as follows: Is any endomorphism $\varphi$ of a differential polynomial algebra $A$ an automorphism if $J(\varphi)$ is invertible over $A[\Delta]$? The study of the Jacobian Conjecture for differential polynomial algebras is interesting first of all in the perspective of constructing a counter example. Partial differential polynomial algebras  have wild automorphisms even in two variables \cite{DNU}.

\section{Differential algebraic dependence}

\hspace*{\parindent}

Elements $f_1,\ldots,f_m$ of $B=k\langle x_1,x_2,\ldots,x_n\rangle$ are called {\em differentially algebraically dependent} if there exists a nonzero element $g$ of a differential polynomial algebra $k\{y_1,\ldots,y_m\}$ in the variables $y_1,\ldots,y_m$ such that $g(f_1,\ldots,f_m)=0$.  Obviously, the elements $f_1,\ldots,f_m$ are differentially algebraically dependent if and only if the elements $f_1^{\theta},\ldots,f_m^{\theta}$, where $\theta\in \Theta$, are algebraically dependent. The set of elements $f_1^{\theta},\ldots,f_m^{\theta}$, where $\theta\in \Theta$, is infinite and we cannot check algebraic dependence of this set by using the algorithms for polynomial algebras.

\begin{theor}\label{t1} Elements $f_1,\ldots,f_m$ of $B=k\langle x_1,x_2,\ldots,x_n\rangle$ are  differentially algebraically dependent if and only if $d(f_1),\ldots,d(f_m)$ are left dependent over $B[\Delta]$.
\end{theor}

\Proof Elements $f_1,\ldots,f_m$ are differentially algebraically dependent if and only if the elements $f_1^{\theta},\ldots,f_m^{\theta}$, where $\theta\in \Theta$, are algebraically dependent. An infinite set of elements is algebraically dependent if and only if it contains a finite algebraically dependent subset. Let $\theta_{1,1},\ldots,\theta_{1,s_1},\theta_{2,1},\ldots,\theta_{m,s_m}$ be elements of $\Theta$ such that $\theta_{i,j}\neq \theta_{i,t}$ if $j\neq t$ and the elements
\bee\label{f1}
f_1^{\theta_{1,1}},\ldots,f_1^{\theta_{1,s_1}},f_2^{\theta_{2,1}},\ldots,f_m^{\theta_{m,s_m}}
\eee
are algebraically dependent.

Let $e_1,\ldots,e_r$ be the minimal subset of $X^{\Theta}$ such that the set of elements (\ref{f1}) is contained in $k[e_1,\ldots,e_r]$. For any $f\in k[e_1,\ldots,e_r]$ let
\bes
\partial(f)=(\frac{\partial f}{\partial e_1},\ldots,\frac{\partial f}{\partial e_r}),
\ees
where $\frac{\partial}{\partial e_i}$ is the usual partial derivation of the polynomial algebras. The elements (\ref{f1}) are algebraically dependent  \cite{SU4} if and only if the elements
\bes
\partial(f_1^{\theta_{1,1}}),\ldots,\partial(f_1^{\theta_{1,s_1}}),
\partial(f_2^{\theta_{2,1}}),\ldots,\partial(f_m^{\theta_{m,s_m}})
\ees
are left dependent over $k[e_1,\ldots,e_r]$. Suppose that
\bee\label{f2}
\sum_{i,j} a_{ij}\partial(f_i^{\theta_{i,j}})=0,
\eee
where $a_{ij}\in k[e_1,\ldots,e_r]$ and at least one $a_{ij}$ is nonzero.

Notice that for any $f\in k[e_1,\ldots,e_r]$ we have
\bes
d(f)=\partial(f)J(e_1,\ldots,e_r).
\ees
Consequently, multiplying (\ref{f2}) by $J(e_1,\ldots,e_r)$ from the right hand side, we get
\bes
\sum_{i,j} a_{ij} d(f_i^{\theta_{i,j}})=0.
\ees
Since $d$ commutes with the elements of $\Delta$, this equality can be written as
\bes
\sum_{i,j} a_{ij}\theta_{i,j} d(f_i)= \sum_{i} (\sum_{j} a_{ij}\theta_{i,j}) d(f_i)=  \sum_{i} b_i d(f_i)=0,
\ees
where $b_i=\sum_{j} a_{ij}\theta_{i,j}$. If $a_{ij}\neq 0$, then $b_i\neq 0$, since $\theta_{i,j}\neq \theta_{i,t}$ if $j\neq t$ and $d(f_1),\ldots,d(f_m)$ are left dependent over $B[\Delta]$. Conversely, every dependence of this form can be uniquely written in the form (\ref{f2}). $\Box$

\begin{lm}\label{l1} Let $k$ be a constructive differential field. Then for any two nonzero elements $a$ and $b$ of the algebra $B[\Delta]$ there can be effectively found two elements $c,d\in B[\Delta]$ such that $ca=db\neq 0$.
\end{lm}
\Proof It is well known that if $R$ is a left Noetherian domain with a derivation $\delta$ then the skew polynomial ring $R[x,\delta]$ is again a left Noetherian domain \cite{GW}.

Let $B_i$ be the subalgebra of $B[\Delta]$ generated by $B$ and $\delta_1,\ldots,\delta_i$, where $0\leq i\leq m$. In particular, $B_0=B$ and $B_m=B[\Delta]$. Notice that $\delta_{i+1}$ can be considered as a derivation of $B_i$ with trivial zero action on $\delta_1,\ldots,\delta_i$. Then $B_{i+1}\cong B_i[\delta_{i+1}]$ is a skew polynomial ring over $B_i$. Consequently, $B[\Delta]$ is a left Noetherian domain.

It is also well known that every left Noetherian domain satisfies the Ore condition, that is, for any nonzero elements $a,b$ there exist $c,d$ such that $ca=db\neq 0$. This proves the existence of elements $c,d$ in the formulation of the lemma. Notice that $a,b$ are left dependent if and only if there exists a nonnegative integer $s$ such that the set of elements $a^{\theta},b^{\theta}$, where $\theta\in \Theta$ and $|\theta|\leq s$, is left dependent over $B$. The latter is algorithmically recognizable since $B$ is a constructive field. Starting from $0$, we can find the minimal $s$ satisfying this condition and a nontrivial equality
\bee\label{f3}
\sum_{|\theta|\leq s} c_{\theta}a^{\theta}=\sum_{|\theta|\leq s} d_{\theta}b^{\theta},
\eee
where $c_{\theta},d_{\theta}\in B$. Let
\bes
c=\sum_{|\theta|\leq s} c_{\theta}\theta,\;\; d=\sum_{|\theta|\leq s} d_{\theta}\theta.
\ees
Then (\ref{f3}) means that $ca=db$. The nontriviality of (\ref{f3}) implies that one of the elements $c$ and $d$ is nonzero. Then  $ca=db\neq 0$ since $B[\Delta]$ is a domain. $\Box$

\begin{lm}\label{l2} Let $k$ be a constructive differential field. Then the left dependency of a finite system of elements of a free left $B[\Delta]$-module $B[\Delta]^s$ is algorithmically recognizable.
\end{lm}

\Proof Let $Q$ be the left classical quotient ring of $B[\Delta]$. Recall that every element of $Q$ can be represented in the form $a^{-1}b$, where $a,b\in B[\Delta]$. We want to show that $Q$ is a constructive skew field. In fact, $(a^{-1}b)^{-1}=b^{-1}a$ and
\bes
a^{-1}bc^{-1}d=a^{-1}c_1^{-1}b_1d=(c_1a)^{-1}(b_1d),
\ees
where $b_1,c_1\in B[\Delta]$ satisfying the condition $bc^{-1}=c_1^{-1}b_1$ can be found effectively by Lemma \ref{l1}. Moreover,
\bes
a^{-1}b+c^{-1}d=t^{-1}(r_1b+r_2d),
\ees
where $t,r_1,r_2\in B[\Delta]$ satisfying condition $t=r_1a=r_2b$ can also be effectively found by Lemma \ref{l1}.

Since $B[\Delta]$ is an Ore ring, it follows that a finite system of elements of the free left $B[\Delta]$-module $B[\Delta]^s$ is left dependent if and only if they are left dependent elements of the free $Q$-module $Q^s$. Using the usual triangulation, we can recognize the left dependency of elements of a finite dimensional vector space over a constructive skew field $Q$. $\Box$

\begin{theor}\label{t2} Let $k$ be a constructive differential field. Then the differential algebraic dependency of a finite system of elements of a free differential field of rational functions over $k$ is algorithmically recognizable.
\end{theor}

\Proof By Theorem \ref{t1}, any elements $f_1,\ldots,f_m$ of $B=k\langle x_1,x_2,\ldots,x_n\rangle$ are  differentially algebraically dependent if and only if $d(f_1),\ldots,d(f_m)$ are left dependent over $B[\Delta]$. The latter is algorithmically recognizable by Lemma \ref{l2}. $\Box$

\section{A representation of free Novikov algebras}

\hspace*{\parindent}

A nonassociative algebra $A=(A, \circ)$ is called a \textit{(left) Novikov algebra} if $A$ satisfies the following identities:
$$(a\circ b)\circ c-a\circ (b\circ c)=(b \circ a)\circ c-b\circ (a\circ c),$$
$$(a\circ b)\circ c=(a\circ c)\circ b,$$
for any $a,b,c\in A $.

Let $k\{x_1,x_2,\ldots,x_n\}$ be the differential polynomial algebra over a field $k$ of a characteristic $0$ in the variables $x_1,x_2,\ldots,x_n$ with one derivation $\delta$. For convenience we denote the derivatives $a^\delta,a^{\delta^2}, a^{\delta^s}$ by $a',a'',a ^{(s)}$, respectively. Put $X=\{x_1,x_2,\ldots,x_n\}$ and by $X^{\delta}$ we denote the set of all symbols of the form $x_i^{(r)}$, where $1\leq i \leq n$, $r \in \mathbb{Z_+}$. For any $x_i^{(r)}, x_j^{(s)}\in X^{\delta}$ we assume that $x_i^{(r)}>x_j^{(s)}$ if $i>j$ or $i=j$, $r>s$. The set $M$ of all differential monomials of the form
$$u = x_{i_1}^{(s_1)}x_{i_2}^{(s_2)} \ldots x_{i_t}^{(s_t)},$$
where $t\geq 0$, $x_{i_j}^{(s_j)} \in X^{\delta}$ for all $1 \leq j \leq t$ and $x_{i_1}^{(s_1)}\geq x_{i_2}^{(s_2)} \geq \ldots \geq x_{i_t}^{(s_t)}$, provides a linear basis for the algebra $ k\{x_1,x_2,\ldots,x_n\}$.

For any $x_i^{(r)} \in X^{\delta}$ let
$$\deg(x_i^{(r)})=1, \; \; \; d(x_i^{(r)})=r,$$
where $1\leq i\leq n$. If $u=a_1\ldots a_s\in M$, where $a_1,\ldots,a_s \in X^\delta$, then let
$$\deg(u)=\deg(a_1)+\ldots+\deg(a_s), \; \; \; d(u)=d(a_1)+\ldots+d(a_s),$$
i.e., $\deg(u)$ is the degree of $u$ with respect to the variables $x_1,\ldots,x_n$ and $d(u)$ is the degree of $u$ with respect to $\delta$. Denote by $\deg_{x_i}(u)$ the degree of the monomial $u$ with respect to $x_i$ ($1\leq i\leq n$).

Let $\rho(u)=\deg(u)-d(u)$ for any differential monomial $u$. It is easy to check that
\bee \label{f4}
\deg(fg)=\deg(f)+\deg (g), \ \
d(fg)=d(f)+d(g), \ \
\rho(fg)=\rho(f)+\rho(g),
\eee
and
\bee \label{f5}
\deg(f^{(r)})=\deg(f), \ \
d(f^{(r)})=d(f)+r, \ \
\rho(f^{(r)})=\rho(f)-r,
\eee
for all $f,g \in k\{x_1,x_2,\ldots,x_n\}$.

The degree function $\rho$ defines a grading
\bes
C=\bigoplus_{i\in \mathbb{Z}}C_i
\ees
of algebra $C=k\{x_1,x_2,\ldots,x_n\}$, where $C_i$ is the $k$-span of differential monomials $u_i$ such that $\rho(u_i)=i$. Every non-zero element $c\in C$ can be uniquely written in the form
\bes
c=c_{i_1}+c_{i_2}+\ldots+c_{i_s}, \ \ i_1<i_2<\ldots<i_s, \ \  0\neq c_{i_j}\in C_{i_j}.
\ees
The element $c_{i_s}$ is called {\em the highest homogeneous part} of  $c$ with respect to the degree function $\rho$ and will be denoted by $\overline{c}$.

\begin{co} \label{c1}
If $f,g \in k\{x_1,\ldots,x_n\}$ are homogeneous elements with respect to degree function $\rho$, then $fg$ and $f^{(r)}$ are also homogeneous with respect to $\rho$.
\end{co}

On the differential polynomial algebra $k\{x_1,\ldots,x_n\}$ we introduce a new operation $\circ$  by assuming
$$f\circ g=fg', \ \ \ f,g\in k\{x_1,x_2,\ldots,x_n\}.$$
Then $k\{x_1,x_2,\ldots,x_n\}$ with the new operation $\circ$ becomes a Novikov algebra. Denote by $N_0 \left\langle x_1,x_2,\ldots,x_n \right\rangle$ the subalgebra of this algebra generated by the elements
$x_1,x_2,\ldots,x_n$.  It is proved in \cite{DzhL} that $N_0\left\langle x_1,x_2,\ldots,x_n \right\rangle$ is a free Novikov algebra in the variables $x_1,x_2,\ldots,x_n$ without identity.

The structure of the space $N_0 \left\langle x_1,x_2,\ldots,x_n \right\rangle$ can be described as follows.

\begin{pr} \cite{DU} \label{p1}
The set of all differential monomials $u\in M$ with the condition $\rho(u)=1$ is a basis of the free Novikov algebra $N_0 \left\langle x_1,\ldots,x_n \right\rangle$.
\end{pr}

 Denote by $N\left\langle x_1,\ldots, x_n\right\rangle=k\oplus N_0\left\langle x_1,\ldots, x_n\right\rangle=N_0\left\langle x_1,\ldots, x_n\right\rangle^{\#}$  the algebra obtained from $N_0 \left\langle x_1,\ldots,x_n \right\rangle$ by formally joining the identity. Then $N \left\langle x_1,\ldots,x_n\right\rangle$ is the free Novikov algebra over $k$ in the variables $x_1,\ldots,x_n$ with identity. This representation of the free Novikov algebra $N \left\langle x_1,\ldots,x_n \right\rangle$ will be used below.

\section{Novikov dependence}

\hspace*{\parindent}

As in the preceeding section, let $k\{x_1,x_2,\ldots,x_n\}$ be the differential polynomial algebra over a field $k$ of a characteristic $0$ in the variables $x_1,x_2,\ldots,x_n$ with one derivation $\delta$.

\begin{lm}\label{l3} Let $f_1,f_2,\ldots,f_p\in k\{x_1,\ldots,x_n\}$ be homogeneous elements with respect to the degree function $\rho$ and $\rho(f_i)=1$ for all $1\leq i\leq p$. If $u\in k\{z_1,\ldots,z_p\}$ is a monomial such that $\rho(u)=s$, then the element $u(f_1,f_2,\ldots,f_p)$ is also homogeneous with respect to the degree function $\rho$ and
$$\rho(u(f_1,f_2,\ldots,f_p))=s.$$
\end{lm}
\Proof If $u=z_i^{(r)}$, then $\rho(u)=1-r=s$. By Corollary \ref{c1} $u(f_i)=f_i^{(r)}$ is homogeneous with respect to $\rho$  and \eqref{f5} gives that
$$\rho(u(f_i))=\rho(f^{(r)}_ i)=\rho(f_i)-r=1-r=s.$$

Let $u=vw$, where $\deg(v),\deg (w)<\deg(u)$. Conducting induction on $\deg(u)$, we can assume that $v(f_1,f_2,\ldots,f_p)$ and $w(f_1,f_2,\ldots,f_p)$ are homogeneous with respect to $\rho$ and $\rho(v(f_1,f_2,\ldots,f_p))=\rho(v)$, $\rho(w(f_1,f_2,\ldots,f_p))=\rho (w)$.
From \eqref{f4} we get
\bes
\rho(u(f_1,f_2,\ldots,f_p))=\rho(v(f_1,f_2,\ldots,f_p)w(f_1,f_2,\ldots,f_p))\\
=\rho(v(f_1,f_2,\ldots,f_p))+\rho(w(f_1,f_2,\ldots,f_p))=\rho(v)+\rho(w)=s.  \ \ \ \Box
\ees

The elements $f_1,f_2,\ldots f_p$ of the Novikov algebra $N \left\langle x_1,x_2,\ldots,x_n \right\rangle$ are called \textit{Novikov dependent} if there exists a nonzero element $h(z_1,z_2,\ldots,z_p)\in N\left\langle z_1,z_2,\ldots,z_p\right\rangle$ such that $h(f_1,f_2,\ldots,f_p)=0$.

\begin{theor} \label{t3}
Let $f_1,f_2,\ldots,f_p$ be elements of the free Novikov algebra $N_0 \left\langle x_1,x_2,\ldots,x_n \right\rangle$ (without identity). The elements $f_1,f_2,\ldots,f_p$ are Novikov dependent if and only if they are differentially algebraically dependent in $k\{x_1,x_2,\ldots,x_n\}$.
\end{theor}
\Proof  Suppose that the elements $f_1,f_2,\ldots,f_p$ of the Novikov algebra $N_0 \left\langle x_1,x_2, \ldots,x_n \right\rangle$ are Novikov dependent. Let   $h(z_1,z_2,\ldots,z_p)\in N \left\langle z_1,z_2,\ldots,z_p \right\rangle$ be a nonzero element such that $h(f_1,f_2,\ldots,f_p)=0$. Since $f_1,f_2,\ldots,f_p$ do not contain constants, we can assume that $h(z_1,z_2,\ldots,z_p)$ also does not contain a constant, i.e., $h(z_1,z_2,\ldots,z_p)\in N_0\left\langle z_1,z_2,\ldots,z_p \right\rangle$. Consequently, $h(z_1,z_2,\ldots,z_p)\in k\{z_1,z_2,\ldots,z_p\}$. Then $h(f_1, f_2,\ldots,f_p)=0$ gives that $f_1,f_2,\ldots,f_p$ are differentially algebraically dependent.

Suppose that $f_1,f_2,\ldots,f_p\in N_0 \left\langle x_1,x_2,\ldots,x_n \right\rangle$ are differentially algebraically dependent. Let $g (z_1,z_2,\ldots,z_p)\in k\{z_1,z_2,\ldots,z_p\}$ be a nonzero element such that $g(f_1,f_2,\ldots,f_p)=0$. We can write $g$ in the form
\bes
g=u_t+u_{t+1}+\ldots+u_s,
\ees
where $\rho(u_i)=i$, $t\leq i\leq s$, and $u_s\neq 0$. Then
$$g(f_1, f_2, \ldots, f_p)=u_t(f_1, f_2, \ldots, f_p)+u_{t+1}(f_1, f_2, \ldots, f_p)+\ldots+u_s(f_1, f_2, \ldots, f_p)=0.$$
By Lemma \ref{l3} we have
$$\overline{g(f_1,f_2,\ldots,f_p)}=u_s(f_1,f_2,\ldots,f_p).$$
It follows that $u_s(f_1,f_2,\ldots,f_p)=0$. If $s>1$, then by differentiating $u_s(f_1,f_2,\ldots,f_p)$ $s-1$ once we get the element $w_1$ such that
$$w_1(f_1,f_2,\ldots,f_p)=(u_s(f_1, f_2, \ldots, f_p))^{(s-1)}=0,\;\;\;\; \rho(w_1)=1.$$
 By Proposition \ref{p1}, the element $w_1$ is an element of the Novikov algebra $N_0\left\langle z_1,z_2,\ldots,z_p \right\rangle$. Then $f_1,f_2,\ldots,f_p$ are Novikov dependent. If $s<1$, then multiplying $u_s(f_1,f_2,\ldots,f_p)$ by $f_1^{\left | s \right|+1}$ we get an element $w_2$ such that
$$w_2(f_1,f_2,\ldots,f_p)=u_s(f_1, f_2, \ldots, f_p) f_1^{\left|s\right|+1}=0,\;\;\;\; \rho(w_2)=1.$$
Therefore, the elements $f_1,f_2,\ldots,f_p$ are again Novikov dependent. $\Box$

\begin{co} \label{c2}
Novikov dependence of a finite system of elements of a free Novikov algebra  over a constructive field $k$ of characteristic zero is algorithmically recognizable.
\end{co}

\Proof Let $f_1,f_2,\ldots,f_p \in N \left\langle x_1,x_2,\ldots,x_n \right\rangle$. We can assume that $f_1,f_2,\ldots,f_p \in N_0 \left\langle x_1,x_2,\ldots,x_n \right\rangle$ since the constants do not affect the dependence of elements. Theorems \ref{t2} and \ref{t3} complete the proof. $\Box$

\begin{co} \label{c3}
Any $n+1$ elements of a free Novikov algebra of rank $n$ over a field of characteristic zero are Novikov dependent.
\end{co}
\Proof It is well known \cite{Kolchin} that any $n+1$ elements of a differential polynomial algebra of rank $n$ are differentially algebraically dependent. $\Box$

\end{document}